\documentclass[12pt, twoside,reqno,a4paper]{amsart}
\usepackage{amsmath}
\usepackage{amsthm}
\usepackage{graphicx}
\usepackage{amssymb}
\usepackage{color}

\newtheorem{thm}{Theorem}[section]

\theoremstyle{definition}

\theoremstyle{remark}

\newcommand{\wbar}{{\overline{w}}}

\begin{document}

\title{A New Proof Of Nakazi Theorem}

\author[Yadav]{Tapesh Yadav}
\address[Yadav]{Department of Mathematics, Indian Institute of Science, Banaglore 560 012, India.}
\email{tapeshyadav1@gmail.com}	
\maketitle

\begin{abstract}
We give a new proof of an important theorem by Nakazi using recent results by Sarason in his seminal paper on agebraic properties of truncated Toeplitz operators.

\end{abstract}

	\section{Introduction}

We follow the notations of Sarason \cite{Sarason}. Let $\mathbb D$ denote the open unit disc, $H^2$ the Hardy space on $\mathbb D$, the Hilbert space of holomorphic functions with square summable Taylor series around $0$. For an inner function $u$, denote by $K^2_u$ the orthocomplement $H^2 \ominus uH^2$ of the shift invariant subspace $uH^2$. All functions in Hardy space are identified with their radial limits defined on the unit circle and hence the space $H^2$ is identified as a subspace of $L^2$ of the unit circle with the normalized arc length measure. So, any input to a function, say $z$ in this paper will denote a point on the unit circle and thus will satisfy $|z|=1$. $P_{H^2}$ will denote the projection from $L^2$ to $H^2$.

	For an $L^\infty$ function $\varphi$ on the unit circle, we write the multiplication operator with symbol $\varphi$ on $L^2$ as $M_\varphi$ where $M_\varphi(f)= \varphi f$, $f \in L^2$ . We will represent the Toeplitz operator with symbol $\varphi$ by $T_\varphi$, i.e., $T_\varphi = P_{H^2}M_\varphi|_{H^2}$. Let $P_c$ denote the projection operator from $H^2$ to the one-dimensional space of constant functions.

	We will be using a tool called conjugation denoted by $C$, which acts on $ \mathcal{K}_u^2$ as
$$C(g)(z) = u(z)\overline{zg(z)} , \;\; g \in   \mathcal{K}_u^2.$$
Conjugation is an antiunitary involution. Conjugation exhibits many other remarkable properties (see \cite[p.~495]{Sarason}). In this paper, we will content ourselves by using the fact that conjugation acts on an element of $\mathcal{K}_u^2$ and gives back an element of $\mathcal{K}_u^2$. $Cf$ will routinely be denoted by $\tilde{f}$.

	If $k_w$ denotes the reproducing kernel $k_w(z) = (1 - z\wbar)^{-1}$ on the Hardy space, then its projection $P_uk_w$ on $\mathcal{K}_u^2$ is denoted by $k_w^u$. Hence $k_0^u = P_uk_0 = P_u1$ as $k_0=1$, and $\tilde{k_0^u} = Ck_0^u$.  

	Before moving forward, it will be worthwhile to mention that whenever $\mathcal{K}_u^2$ is non trivial (i.e., is a proper non zero subspace of $H_2$), then $k_0^u \neq 0$ and consequently $\tilde{k_0^u} \neq 0$. Indeed, as $k_0^u(z) = 1- \overline{u(0)}u(z)$, (see \cite[p.~494]{Sarason}), if $k_0^u$ is constant, then $u$ has to be constant, which can never give rise to a non trivial $\mathcal{K}_u^2$. It is clear that if the inner function $u$ is not constant, then $\mathcal{K}_u^2$ is non trivial and hence $\tilde{k_0^u} \neq 0$.\\

A new proof of the following theorem will be the main content of this paper. Our technique relies heavily on results from a recent paper by Sarason (see \cite{Sarason}).  We state here the theorem by Takahiko Nakazi (see \cite{Nakazi}) which he proved in much more generality.

	\begin{thm}[Nakazi Theorem] \label{nakathm}
		Let $\varphi$ be an $L^\infty$ function on the unit circle. Decompose $\varphi$ as $\bar{f}+zg$, where $f,g \in H^2 $. Then, for any non constant inner function $u$, $uH^2$ is invariant under $T_\varphi$ if and only if $f$ is a constant.
	\end{thm}
	
	Recall that by Beurling's theorem, any closed subspace of $H^2$ that is invariant under the forward shift is of the form $uH^2$ for some inner function $u$. Thus, a non trivial invariant subspace of the shift is invariant under a Toeplitz operator if and only if the symbol of the Toeplitz operator is analytic.

	\section{Proof of Nakazi Theorem}

	\begin{proof}
		Pick an element from $uH^2$ say $uh$, where $h \in H^2$. We will first show that if $f(z) = c$, where $c$ is a constant, then, $uH^2$ is invariant under $T_\varphi$. So, we compute:
		\begin{equation} \label{eq1}
			T_\varphi(uh) = T_{\bar{f}+zg}(uh) = T_{\bar{c}}(uh) + T_{zg}(uh).
		\end{equation}
		Clearly, the second term in the last expression i.e., $T_{zg}(uh) $ lies in $uH^2$ as $T_{zg}(uh) = P_{H^2}(u(zgh))= u(zgh) \in uH^2$. Similarly, the first term in (\ref{eq1}) i.e., $T_{\bar{c}}(uh) = P_{H^2}(\bar{c}uh) = \bar{c}uh \in uH^2$. Thus, $T_{\bar{c}}(uh) + T_{zc}(uh) \in uH^2$ and hence  $T_\varphi(uh) \in uH^2$. This shows that $uH^2$ is invariant under $T_\varphi$.\\
		
		Now, we will show that  $uH^2$ is not invariant under $T_\varphi$ whenever $f$ is non constant and vanishes at $0$, i.e., if $f$ is orthogonal to the constant functions. Decompose such an $f$ as $f_1+f_2+f_3$ according to the decomposition  $H^2 = \mathcal{K}_u^2 \oplus uH^2 = span\{k_0^u\} \oplus (\mathcal{K}_u^2 \ominus span\{k_0^u\} ) \oplus uH^2$ of the space, i.e., $f_1 \in span\{k_0^u\}$, $f_2 \in \mathcal{K}_u^2 \ominus span\{k_0^u\}$ and $f_3 \in uH^2$. Thus $f_1 = ck_0^u, c \in \mathbb{C}$ and $f_3 = ul, l \in H^2$. Our strategy to prove the claim will be as follows. For any $h \in H^\infty$, $ T_\varphi(uh) = P_{H^2}(\bar{f}(uh)) + P_{H^2}(zg(uh))$. Since $P_{H^2}((zg)uh) = u(zgh) \in uH^2$. Thus, if we show that $P_{H^2}(\bar{f}(uh)) \not \in uH^2$ for some $h \in H^\infty$, then $T_\varphi(uh) \not \in uH^2$, hence proving our claim.\\ 

		\textbf{Case 1:} $\mathbf{f_2 \neq 0}$ \\
		Note that, for any $h \in H^\infty$, we have
		\begin{align*} \label{9a}
		P_{H^2}(\bar{f}(uh)) =  P_{H^2}((\overline{f_1+f_2+f_3})(uh))
		&= P_{H^2}(\overline{ck_0^u}uh + \bar{f_2}uh + \overline{ul}uh)\\ &=  P_{H^2}(\overline{c} \tilde{k_0^u}zh + \tilde{f_2}zh + \overline{l}h).
		\end{align*}
		The last equality follows from the definition of conjugation and the fact that an inner function is 1 a.e. Observe that $\overline{c} \tilde{k_0^u}zh$ and $ \tilde{f_2}zh$ lie in $H^2$. So,
		\begin{equation}\label{mainq}
		    P_{H^2}(\bar{f}(uh)) =  \overline{c} \tilde{k_0^u}zh + \tilde{f_2}zh + P_{H^2}(\overline{l}h)
		\end{equation} 
		If we take $h=1$, we get
	    \begin{equation*}
		    P_{H^2}(\bar{f}u) =  \overline{c} \tilde{k_0^u}z + \tilde{f_2}z + P_{H^2}(\overline{l})
		\end{equation*} 
		Since $\overline{l}$ is coanalytic, $P_{H^2}(\overline{l}) = c_1$, a constant. Also, note that $\tilde{f_2}z \in \mathcal{K}_u^2$ as $\tilde{f_2} \in \mathcal{K}_u^2 \ominus   \text{span}\{ \tilde{k_0^u} \}$ (see remark in \cite[p.~512]{Sarason}). Hence we have
		\begin{align*}
		\langle P_{H^2}(\bar{f}u), \tilde{f_2}z \rangle &= \langle \overline{c} \tilde{k_0^u}z + \tilde{f_2}z + P_{H^2}(\overline{l}), \tilde{f_2}z \rangle = \langle \overline{c} \tilde{k_0^u}z, \tilde{f_2}z \rangle + \langle\tilde{f_2}z, \tilde{f_2}z \rangle + \langle  P_{H^2}(\overline{l}), \tilde{f_2}z \rangle.
		\end{align*}
		Note that, $ \langle c_1, \tilde{f_2}z \rangle = 0$ and, $ \langle \overline{c} \tilde{k_0^u}z, \tilde{f_2}z \rangle = \langle \overline{c} \tilde{k_0^u}, \bar{z}\tilde{f_2}z \rangle = \langle \overline{c} \tilde{k_0^u}, \tilde{f_2} \rangle = 0$ where the first equality follows from the fact that $M_z^\star = M_{\bar{z}}$. The second equality follows from the fact that $z\bar{z} = 1$ and the third equality follows as $\tilde{f_2} \in \mathcal{K}_u^2 \ominus   \text{span}\{ \tilde{k_0^u} \}$. So, we get,
		$$ \langle	P_{H^2}(\bar{f}u), \tilde{f_2}z \rangle = {|| \tilde{f_2}z ||}^2 $$
		which is a non-zero because of the assumption that $f_2 \neq 0$. Thus $P_{H^2}(\bar{f}u) \not\in uH^2 $ as $\tilde{f_2}z \in {uH^2}^{\perp}$ and hence the claim.\\

		\textbf{Case 2:} $\mathbf{f_2 = 0}$ \\
		Now we shall assume that $f_2 = 0$ i.e., $f = f_1 + f_3 = ck_0^u + ul$ with all the symbols same as in Case 1. \\ \\
		If $c=0$, then equation (\ref{mainq}) gives that $ P_{H^2}(\bar{f}(uh)) =  P_{H^2}(\overline{l}h)$. Since $l\in H^2$, $l$ admits a series expansion of the form $l= \sum_{n=N}^{\infty} a_n z^n$, where $a_N$ is non zero and $N \geq 0$. So, $\overline{l} = \sum_{n= N}^{\infty} \overline{a_n} z^{-n}$ and therefore $\overline{l} z^N = \sum_{n= N}^{\infty} \overline{a_n} z^{-n+N} = \sum_{m=0}^{\infty} \overline{a_{m+N}} z^{-m}$. If we take $h_N(z)=z^N$, then $P_{H^2}(\overline{l}h_N) = \overline{a_N} \neq 0$. Hence, $$ \langle P_{H^2}(\bar{f}(uh_N)), k_0^u \rangle = \langle \overline{a_N}, k_0^u \rangle = \overline{a_N} \langle 1, k_0^u \rangle = \overline{a_N} {|| k_0^u||}^2 \neq 0$$
		Since $k_0^u \in \mathcal{K}_u^2 = {uH^2}^{\perp}$, we get that $P_{H^2}(\bar{f}(uh_N)) \not\in uH^2 $  which proves the claim that $uH^2$ is not invariant under $T_\varphi$.\\ \\
		If $c\neq 0$, we can assume without loss of generality that $f= k_0^u + ul$. Note that for $f$ to belong to $H^2 \ominus   \text{span} \{ 1 \}$, we need that $uP_c(l) \neq 1- k_0^u$, which gives $ P_c(l) \neq \overline{u(0)}$, which follows from the fact that $ k_0^u(z) = 1 - \overline{u(0)}u(z)$ \cite[p.~494]{Sarason}. From equation (\ref{mainq}), we have
		\begin{equation} \label{eq2}
		P_{H^2}(\bar{f}u) = P_{H^2}(\bar{k_0^u}u + \overline{ul}u) = \tilde{k_0^u}z + P_{H^2}(\bar{l})
		\end{equation}
		Now, we shall make an important observation that $\tilde{k_0^z} = T_z^\star(u)$ \cite[p.~495]{Sarason}. And $T_z^\star(u) = \frac{u(z) - u(0)}{z}$, which is direct from definition of backward shift $T_z^\star$ \cite[p.~492]{Sarason}. Thus, we get that $ \tilde{k_0^u}z = u(z) - u(0)$. So, (\ref{eq2}) becomes,
		$$ P_{H^2}(\bar{f}u) = (u - u(0)) + P_{H^2}\bar{l} = (u - u(0)) + P_{c}\bar{l}$$
		The last equality follows from the fact that $\bar{l}$ is co analytic. Call $c_2 = P_{c}\bar{l}$. Clearly, $c_2 \neq u(0)$ as is shown above. Now, taking inner product with $k_0^u \in \mathcal{K}_u^2$,
		$$ \langle P_{H^2}(\bar{f}u), k_0^u \rangle = \langle (u - u(0)) + c_2, k_0^u \rangle = \langle u, k_0^u \rangle  + \langle c_2 - u(0), k_0^u \rangle = (c_2 - u(0)) {||k_0^u||}^2 \neq 0$$
		as $(c_2 - u(0)) \neq 0$ from last paragraph. Hence $P_{H^2}(\bar{f}u) \not\in uH^2 $ which proves the claim that $uH^2$ is not invariant under $T_\varphi$ in this case as well. \\ \\
		From the above two cases, it follows that for all $f$ which are in the orthocomplement of the constant function $1$, $uH^2$ is not invariant under $T_\varphi$. Since $f \in H^2$ can be written as $f = f_4 + f_5$ where $f_4 \in   \text{span}\{ 1 \}$ and $f_5 \in H^2 \ominus   \text{span}\{ 1 \}$, we have $\varphi = \bar{f} + zg = (\overline{f_4 + f_5}) + zg= \bar{f_4} +  (\bar{f_5} + zg)$. Hence, we get $T_\varphi =  T_{\bar{f_4}}  + T_{\bar{f_5} + zg}$. Note that $T_{\bar{f_4}}$ leaves $uH^2$ invariant. And, if ${f_5} \neq 0$, then $T_{\bar{f_5}+zg}$ does not leave $uH^2$ invariant as we have already shown using Cases 1 and 2. Therefore, $uH^2$ is invariant under $T_\varphi$ if and only if $f$ is a constant, i.e., $f_5 = 0$. 	\end{proof}


\begin{thebibliography}{99}




\bibitem{Nakazi}
{\sc T.~Nakazi}, {\em {Invariant Subspaces Of Toeplitz Operators And Uniform
  Algebras}}, Bull. Belg. Math. Soc. Simon Stevin, 15 (2008), pp.~1--8.

\bibitem{Sarason}
{\sc D.~Sarason}, {\em Algebraic properties of truncated toeplitz operators},
  Oper. Matrices, 1 (2007), pp.~491--526.

\end{thebibliography}
\end{document}